\documentstyle[12pt]{article}
\begin{document}

\title{Fractional spin through \\quantum strange superalgebra $\tilde{P}_{Q}(n)$.}
\author{\bf { Mostafa Mansour \thanks{email: mostafa.mansour.fpb@gmail.com}
  } \\
  D\'epartement de Physique, \\
  Facult\'e Polydisciplinaire de Beni Mellal.\\ Universit\'e Sultan Moulay Slimane. \\Beni Mellal. Morocco.\\
 }
\date{}
\maketitle

\begin{abstract}
The purposes of this paper is to investigate the properties of the quantum extended
strange superalgebra $\tilde{P}_{Q}(n)$ when his deformation parameter $Q$ goes to a root
of unity.
\end{abstract}

\newpage\

\section{Introduction}
In recent years, much interest has been made in the study of the
Lie superalgebras \cite{r01,r02,r03,r04,r05,r3}.
These structures can be obtained through consistent realization involving deformed
Bose and Fermi operators \cite{r1,r2}.\\

In another vein, the geometric interpretation of fractional supersymmetry
 have  gaining increased attention, particularly in the works \cite{r11,r12,r13,r14,r15, r151,r17}, where the authors
  show that the one-dimensional superspace is equivalent to the braided line
  when the deformation parameter goes to a root of unity $Q \rightarrow q$, and the braided line  is generated by a generalized
odd variable and a (classical)ordinary even variable. In the work \cite{r18}, R. S. Dunne, using $Q$-oscillator realization,
proved that the $U_Q(sl(2))$ is similar  to a direct product of the finite
classical algebra $U(sl(2))$ and the q-deformed one $U_q(sl(2))$ (where q is a root of unity).\\

Since there exist $Q$-oscillator
realization of all deformed algebras and deformed super algebras $U_Q(g),$ it is convenable
to explore the  splitting of these (super) algebras when $%
Q\rightarrow q.$ In this context, the property of splitting of  some  particular quantum (Super)-algebras was examined in \cite{r19}. The decomposition of the quantum (super) Virasoro algebras is described in \cite{m2}. The case of quantum affine algebras with vanishing central charge is developed in \cite{m3} and the case quantum algebras $A_{n},  B_{n},  C_{n} $ and $D_{n}$ and quantum superalgebra $A(m,n),  B(m,n),  C(n+1) $ and $D(n,m)$ in the
  $ Q \to q ~limit$ is investigated in \cite{M1}.\\

The Lie superalgebras of classical type are one  of the two
following classes\cite{kac}: basic Lie superalgebras or strange ones. The basic Lie superalgebras
 have proprieties like as simple Lie algebras. They have an invariant non-degenerate bilinear form, but strange Lie superalgebras $P(n)$ and $Q(n)$ haven’t.\\

The strange Lie
superalgebra $P(n)$ have attracted a particularly attention.
In \cite{S1}, Dynkin-like diagrams
of the strange superalgebra
$P(n)$ was examined by Frappat, Sciarrino and Sorba. In \cite{S6},  polynomial representations of strange  Lie superalgebras
are investigated. The oscillator realization of the strange superalgebras $P(n)$  has been constructed by Frappat, Sciarrino and Sorba in \cite{FS}. A deformation   $U_{Q}(\tilde{P}(n)) = \tilde{P}_{Q}(n)$ of the extended non-contragredient (strange) superalgebra $\tilde{P}(n)$ is given in \cite{FSS}.\\

 The purpose of this  paper is to expore the property
  of decomposition  of the quantum extended non-contragredient (strange) superalgebra $U_{Q}(\tilde{P}(n)) = \tilde{P}_{Q}(n)$ in the $ Q \to q ~limit$. In next section (section 2) we review
  some  results  concerning k-fermions, decomposition property of $Q$-boson
  in the $ Q \to q ~ limit$ and the equivalence between $Q$-deformed fermions and  classical ones. Using these results detailed in \cite{m2,m3}, we analyse the $ Q \to q  ~ limit$ of the quantum $ U_{Q}(sl(n))$ algebra of the $sl(n)$ algebra (The bosonic part of $ P(n)$)  (section 3) and the quantum extended non-contragredient (strange) superalgebra of $U_{Q}(\tilde{P}(n))$ (section 4). In the last section (section 5) we shall give some concluding remarks

\section{Preliminaries}
In this section we recall some basic facts about $k$-fermions\cite{r181}, decomposition property of $Q$-boson
  in the $ Q \to q ~ limit$ and the equivalence between  $Q$-deformed fermions and ordinary ones  ( see \cite{m2,m3} for more details).\\\\
Let us began by giving the definition of the Q-bosonic algebra noted ($\Xi^{i} _Q)$, generated by a number operator $N_{A_{i}}$, a creation operator $A_{i}^{+}$ and an annihilation operator $A_{i}^{-}$,
satisfying the  relations:

\begin{equation}
\begin{array}{c}
A_{i}^{-}A_{i}^{+}-Q^{\pm}A_{i}^{+}A_{i}^{-}=Q^{\mp N_{A_{i}}}\\
Q^{N_{A_{i}}}A_{i}^{\pm}Q^{-N_{A_{i}}}=Q^{\pm}A_{i}^{+}\\
Q^{N_{A_{i}}}Q^{-N_{A_{i}}}=Q^{-N_{A_{i}}}Q^{^{+}N_{A_{i}}}=1
\end{array}
\end{equation}
then if we put the following operators as given in \cite{r18}:

\begin{equation}\label{e2}
a_{i}^{-}=\lim_{Q\rightarrow q}\frac{Q^{\pm \frac{kN_{a_{i}}}2}}{([k]!)^{\frac 12}}%
(A_{i}^{-})^k,\,a_{i}^{+}=\lim_{Q\rightarrow q}\frac{(A_{i}^{+})^kQ^{^{\pm }\frac{kN_{A_{i}}}2}%
}{([k]!)^{\frac 12}},
\end{equation}
we can easily show that the above  operators (\ref{e2}) gratifies the relations of an ordinary boson algebra noted $\Xi^{i} _0$, defined  by:

\begin{equation}
\begin{array}{c}
\lbrack a_{i}^{-},a_{i}^{+}]=1.\\
\lbrack N_{a_{i}},a_{i}^{\pm}]=\pm a_{i}^{\pm}
\end{array}
\end{equation}
The number operator of this new algebra is defined  by $%
N_{a_{i}}=a_{i}^{+}a_{i}^{-}$.\\\\
In order  to discuss the splitting
of $Q$-deformed boson in the $limit$~ $Q\rightarrow q$, we introduce the
new operators:
\begin{equation}
\begin{array}{c}
\chi_{i}^{-}=A_{i}^{-}q^{-\frac{kN_{a_{i}}}2}\\
\chi_{i}^{+}=A_{i}^{+}q^{-\frac{kN_{a_{i}}}2}\\
N_{\chi_{i}}=N_{A_{i}}-kN_{a_{i}},
\end{array}
\end{equation}
which satisfies the relations of  a $k$-fermionic algebra noted ($\Sigma^{i} _q)$ defined by
\begin{equation}
\begin{array}{c}
\lbrack \chi_{i}^{+},\chi_{i}^{-}]_{q^{-1}}=q^{N_{\chi_{i}}}\\
\lbrack \chi^{-},\chi^{+}]_q=q^{-N_{\chi_{i}}}\\
\lbrack N_{\chi_{i}},\chi_{i}^{\pm }]=\pm \chi_{i}^{\pm }.
\end{array}
\end{equation}
where the deformation parameter $q=e^{\frac{2i\pi }r},\; r\in N-\{0,1\}$,is a root of unity.\\\\
It straightforward  to check that the two algebras generated respectively by the set of operators
$\{a_{i}^{+},a_{i}^{-},N_{a_{i}}\}$ and $\{\chi_{i}^{+},\chi_{i}^{-},N_{\chi_{i}}\}$ are mutually commutative. We
conclude that in the $limit$ ~$Q\rightarrow q$ , the $Q$-deformed bosonic
algebra oscillator $\Xi^{i}_Q$ decomposes into two independent algebras, an ordinary
boson algebra $\Xi^{i}_0$ and $k$-fermionic algebra $\Sigma^{i} _q$; formally one can write:
\[
\lim_{Q\rightarrow q}\Xi^{i}_Q\equiv \Xi^{i}_0\otimes \Sigma^{i} _q,
\]
We define also the $Q-$deformed fermionic algebra noted $\Omega _Q$ generated by the generators $\Phi_{i}^-,\Phi_{i}^+$ and $Q^{M_{\Phi_{i}}},Q^{-M_{\Phi_{i}}} $ satisfying the following relations
 \begin{equation}
 \begin{array}{c}
 Q^{M_{\Phi_{i}}} Q^{-M_{\Phi_{i}}} = Q^{-M_{\Phi_{i}}} Q^{M_{\Phi_{i}}} = 1\\
  Q^{M_{\Phi_{i}}} Q^{M_{\Phi_{j}}} = Q^{M_{\Phi_{j}}} Q^{M_{\Phi_{i}}}\\
 Q^{M_{\Phi_{i}}}\Phi_{i}^\pm Q^{-M_{\Phi_{i}}} = Q^{\pm\delta_{ij}} \Phi_{j}^\pm \\
 \Phi_{i}^- \Phi_{i}^{+} + Q^{\pm} \Phi_{i}^+ \Phi_{i}^- = Q ^{\pm M_{i}} \\
 \{\Phi_{i}^\pm, \Phi_{i}^\pm\}=0;   \\
 (\Phi^+)^2=0,(\Phi^-)^2=0\\
  \{\Phi_{i}^+, \Phi_{j}^-\}=0 ~~for ~i \neq j
 \end{array}
 \end{equation}
 then if put the new generators
 \begin{equation} \label{eq5}
 \begin{array}{c}
 \phi_{i}^-~=~Q^{-M_{\phi_{i}} \over 2}\Phi_{i}^-\\
 \phi_{i}^+~=~\Phi_{i}^+Q^{-M{\phi_{i}} \over 2}
 \end{array}
 \end{equation}
 we see that the $Q-$deformed fermion reproduces an ordinary fermion  algebra defined by the new operators  (\ref{eq5}) and the following relations
 \begin{equation}
 \begin{array}{c}
 (\phi_{i}^-)^2~=~0\\
 (\phi_{j}^+)^2~=~0\\
 \{ \phi_{i}^-,\phi_{j}^+\}~=~\delta_{ij}
 \end{array}
 \end{equation}
\section{The quantum algebra $U_{Q}(sl(n))$}
Let $C = [ a_{ij}] ( 1 \leq i,j \leq n)$ be a symmetrisable generalized Cartan matrix and let $d_{i} ( 1 \leq i \leq n )$ be the non integers such that
$ d_{i} a_{ij} = a_{ij} d_{i} $. Let $ Q \neq 0 $ be a complex number. For $Q$  generic the quantum enveloping algebra corresponding to $ [ a_{ij}]$ is a Hopf algebra with 1 and generators $ \{ E_{i}, F_{i}, K_{i}^{+} = Q^{+d_{i}H_{i}}, K_{i}^{-} = Q^{-d_{i}H_{i}}; 1 \leq i \leq n\}$ satisfying the following relations :
\begin{equation}
 \begin{array}{c}
 [E_{i}, F_{j}] = \delta_{ij} \frac{K_{i}- K_{i}^{-1}}{Q^{d_{i}} - Q^{-d_{i}}}\\
 K_{i}E_{j} K_{i}^{-1}= Q^{d_{i}a_{ij}}E_{j} \\
  K_{i}F_{j} K_{i}^{-1}= Q^{-d_{i}a_{ij}} F_{j}\\
 K_{i}^{-} K_{i}^{+} = K_{i}^{+} K_{i}^{-}; K_{i} K_{j} = K_{j} K_{i} \\
 \end{array}
 \end{equation}
with Serre relations,
\begin{equation}
\begin{array}{c}
\sum_{0\le t \le n}~(-1)^t~[^n_t]_Q ~(E_i)^t~ E_j~(E_i)^{n-t}=0,\hskip 1cm
i{\not=}j\\
\sum_{0\le t \le n} ~(-1)^t~[^n_t]_Q ~(F_i)^t~ F_j~(F_i)^{n-t}=0,\hskip 1cm
i{\not=}j
\end{array}
\end{equation}
where $a_{ij}$ is Cartan matrix,
$n=1-a_{ij}$ and
\begin{equation}
[^n_t]_Q={\frac{[n]_Q!}{[t]_Q![n-t]_Q!}},\hskip
1.5cm{[t]_Q!=[t]_Q[t-1]_Q...[1]_Q},
\hskip 1.5cm{[t]_Q=\frac{Q^t-Q^{-t}}{Q-Q^{-1}}}
\end{equation}
The Cartan matrices of classical type $A_{n}, B_{n}, C_{n}$ and $ D_{n}$ and the corresponding non zero integers are given in \cite{CP}. A discussion of the Q-boson and Q-fermion representation was given by Hayashi \cite{HA}.\\
We focus here on the algebra $U_{Q}(sl(n))$ (where $sl(n)$ is the bosonic part of $ P(n+1)$)and we assume that the deformation parameter Q is generic. The explicit expressions for corresponding generators as linear and bilinear in Q-deformed bosonic  operators are given by:

\begin{equation}\label{e4}
 \begin{array}{c}
E_{i}= A_{i}^{-} A_{i+1}^{+}\\
F_{i}= A_{i}^{+} A_{i+1}^{-}\\
H_{i} = - N_{A_{i}} + N_{A_{i+1}}
\end{array}
 \end{equation}\\
Now we can explore the $limit ~ Q \rightarrow q $ of the quantum algebra $U_{Q}(sl(n))$. The key tool to discuss this $limit$ is the $Q$-bosonic decomposition presented in detail in \cite{m2,m3}
 when the deformation parameter $Q$ goes to a root of unity $q$. So, the $n$ $Q$-bosons reproduce $n$ ordinary bosons and $n$ $k$-fermions $\{ \chi_{i}^{+}; \chi_{i}^{-}; N_{\chi_{i}} \}$ with $(1 \leq i \leq n- 1$). The $n$ classical bosons are defined by
\begin{equation}\label{e40}
a_{i}^{-}=\lim_{Q\rightarrow q}\frac{Q^{\pm \frac{kN_{A_{i}}}2}}{([k]!)^{\frac 12}}%
(A_{i}^{-})^k,\,a_{i}^{+}=\lim_{Q\rightarrow q}\frac{(A_{i}^{+})^kQ^{^{\pm }\frac{kN_{A_{i}}}2}%
}{([k]!)^{\frac 12}},
\end{equation}
where their number operators are given by $N_{a_{i}} = a^{+}_{i}a^{-}_{i}, for i = 1,2,...,n$. \\
Then, using operators (\ref{e40}), we can construct the classical $U(sl(n))$ algebra(with $1 \leq i \leq n- 1$):
\begin{equation}
 \begin{array}{c}
e_{i}= a_{i}^{-} a_{i+1}^{+}\\
f_{i}= a_{i}^{+} a_{i+1}^{-}\\
h_{i} = - N_{a_{i}} + N_{a_{i+1}}\\
\end{array}
 \end{equation}
From the remaining  operators $\{ \chi_{i}^{+}; \chi_{i}^{-}; N_{\chi_{i}} \}$ with $(1 \leq i \leq n- 1$)we construct the new generators
\begin{equation}
 \begin{array}{c}
E_{i}= \chi_{i}^{-} \chi_{i+1}^{+}\\
F_{i}= \chi_{i}^{+} \chi_{i+1}^{-}\\
H_{i} = - N_{\chi_{i}} + N_{\chi_{i+1}}
\end{array}
 \end{equation}
which realize the $U_{q}(sl(n))$ algebra; where $U_{q}(sl(n))$ is the same version of $U_{Q}(sl(n))$ obtained by simply setting $ Q= q$ rather than by taking the $limit$ as above. The elements of $U_{q}(sl(n))$ and $U(sl(n))$ algebras are mutually commutative.\\
Then, in the $Q \rightarrow q  ~limit$, the quantum algebra $U_{Q}(sl(n))$ is a direct product of the form
\begin{equation}
lim_{Q \rightarrow q} U_{Q}(sl(n)) = U_{q}(sl(n))\otimes U(sl(n))
\end{equation}
Note that, the above direct product  $lim_{Q \rightarrow q} U_{Q}(g) = U_{q}(g)\otimes U(g)$  valid for quantum algebras does not appear of a quantum superalgerbra  $U_{Q}(\mathbf{g})$. In fact the explicit expressions of the generators of the quantum superalgebra $U_{Q}(\mathbf{g})$ are presented as linears and bilinears in Q-deformed bosonic and fermionic oscillator operators, and using the fact that the $Q$-deformed bosonic
operators $\{A^{+}_{i},A^{-}_{i},N_{A_{i}} \}$ with $( 1 \leq i \leq n )$  decomposes into two independent oscillators algebras:  classical
bosons $%
\{a^{+}_{i},a^{-}_{i},N_{a_{i}}\}$ and $k$-fermion operators $\{\chi^{+}_{i},\chi^{-}_{i},N_{\chi_{i}}\}$, and $Q$-fermions become $q$-fermions which are object equivalent to
conventional fermions $\{\phi^{+}_{i},\phi^{-}_{i},M_{\phi_{i}}\}.$
 Then from the classical bosons $\{a_{i}^{-}, a_{i}^{+}, N_{a_{i}} \}$ and classical fermions $\{\phi_{i}^{-}, \phi_{i}^{+}, M_{\phi_{i}} \}$, one can realize
  the nondeformed  superalgebra $U(\mathbf{g})$ but from the remaining operators $\{\chi_{i}^{-}, \chi_{i}^{-}, N_{\chi_{i}} \}$ we construct the generators of a different quantum $q$-algebra ( see \cite{M1} for more details).

\section{The quantum extended strange superalgebra $\tilde{P}_{Q}(n)$}

Let  $\emph{G} =  \emph{G}_{0}\oplus \emph{G}_{1}$ be a $ \texttt{Z}_2$-graded vector space  with $\dim \emph{G}_{0} = j$ and $\dim \emph{G}_{1} = i$. Then there exists a natural superalgebra structure
on the algebra $End\, \emph{G}$ defined by:
\[
End\, \emph{G} = End_{0} \emph{G} \oplus End_{1} \emph{G} ~~\quad where
End_k \emph{G} = \{ \phi \in End\,\emph{G} ~ | ~ \phi ( \emph{G}_l) \subset
\emph{G}_{k+l} \}
\]
The superalgebra
$End\, \emph{G}$ supplied with the Lie superbracket is the Lie superalgebra noted $\ell(i,j).$\\
The elements $M$ of $\ell(i,j)$  have the form
\[
M = \left(\begin{array}{cc} N & Q \cr R & S \end{array}\right)
\]
where $N$ and $S$ are $gl(j)$ and $gl(i)$ matrices, $Q$ and $R$ are
$j \times i$ and $i \times j$ rectangular matrices.
\\\\
The superalgebra of matrices $M \in \ell(n,n)$ satisfying the following equalities
\[
N^t = -S \,, \qquad Q^t = Q \,, \qquad R^t = -R \,, \qquad \mbox{tr}(N) = 0
\]
is nothing other than the (non contragredient) strange superalgebra $P(n)$.\\\\
An oscillator realization of the generators of $P(n)$ given in \cite{r02,FS}. In the Chevalley basis, the (non contragredient) strange superalgebra $P(n)$ is spanned by the generators $ \{ X_{i}, Y_{i}, T_{i}, X_{n} \}$ with $( 1 \leq i \leq n-1 )$ satisfying the following commutation relations:

\begin{equation}
\begin{array}{c}
\left[X_{i}, Y_{j}\right] = \delta_{ij} T_{i}\\
\left[X_{i},T_{j}\right] = - a_{ij}X_{i}\\
\left[Y_{i},T_{j}\right] = a_{ij}Y_{i}\\
\left[T_{i},X_{n}\right] = a_{in}X_{n}\\
\left[T_{i},T_{j}\right] = 0
\end{array}
\end{equation}
where $ (a_{ij})_{1\leq i,j\leq n}$ is the Cartan matrix of $ su(n)$ and $ a_{in}=0$ for $ 1 \leq i \leq n-2$ and $ a_{n-1,n} = -2$.\\\\
Let us to precise that the notions of  Dynkin diagram and Cartan matrix
are not  determined for the non-contragredient  Lie superalgebra $P(n)$. However, if we extend the superalgebra of $P(n)$ by
an appropriate diagonal matrices, one can obtain a non-null bilinear form on the Cartan
subalgebra of this extension and therefore  get in this case a generalized form of  the notions of Cartan
matrix and Dynkin diagram \cite{r02}.\\\\
The extended strange superalgebra $\tilde{P}(n)$ is defined in this basis by
$3(n - 1)$ bosonic generators $ X_{i}, Y_{i}, T_{i}$, with $ i= 1, . . . , n - 1,$ and a fermionic generator $X_{n}$ and a diagonal generator $D$ such that
\begin{equation}
\begin{array}{c}
\left[X_{i}, Y_{j}\right] = \delta_{ij} T_{i}\\
\left[X_{i},T_{j}\right] = - a_{ij}X_{i}\\
\left[Y_{i},T_{j}\right] =  a_{ij}Y_{i};\\
\left[T_{i},X_{n}\right] =  a_{in}X_{n};\\
\left[T_{i},T_{j}\right] = 0;\\
\left[D,X_{i}\right]= [D,Y_{i}]= [D,T_{i}] = 0 \\
 \left[D,X_{n}\right]= X_{n}
\end{array}
\end{equation}
\medskip
The Cartan matrix of the extended strange superalgebra $(a_{ik})_{\tilde{P}(n)}$ of $\tilde{P}(n)$ with $ 1 \leq i \leq n-1$ and $ 1 \leq k \leq n$ is given by:
$$
  (a_{ik})_{\tilde{P}(n)} = \left(\begin{array}{rrrrrrrrrrrr}
 2 & -1 & 0 & \cdots & \cdots &&&  0& 0 \\
 -1 & 2 & -1  &&&&& 0& 0 \\
 0 & -1 & \ddots & \ddots & \ddots &&& \vdots& \vdots \\
 \vdots & 0 & \ddots & \ddots &   \ddots &&&0& 0\\
 \vdots  &&&&&&&\\
 &&&&&& 2 &-1& 0\\
  0&0& &&\cdots&0&-1&2& -2\\

 \end{array}\right)
 $$\\
The Cartan matrix is well defined to get the Serre relations for
the extended non-contragredient  Lie superalgebra $\tilde{P}(n)$ in the quantum case and permits to define a quantum superalgebra structure on the $Q$-deformed version of the  extended non-contragredient (strange) superalgebra $\tilde{P}(n)$.\\
For $Q$ generic, a quantum deformation   $U_{Q}(\tilde{P}(n)) = \tilde{P}_Q(n)$ of the extended non-contragredient (strange) superalgebra $\tilde{P}(n)$ is proposed in \cite{FSS} as follows:
\begin{equation}
\begin{array}{c}
[ \hat{X}_{i}, \hat{Y}_{j} ] = \delta_{ij} [\hat{T}_{i}]_{Q}\\
\left[\hat{T}_{i},\hat{X}_{n}\right] =  a_{in}\hat{X}_{i}\\
\left[\hat{Y}_{i},\hat{T}_{j}\right] =  a_{ij}\hat{Y}_{i}\\
\left[ \hat{T}_{i},\hat{T}_{j} \right] = 0\\
\left[ \hat{D},\hat{X}_{i} \right]= \left[\hat{D},\hat{Y}_{i}\right]= \left[\hat{D},\hat{T}_{i}\right] = 0 \\
\left[ \hat{D},\hat{X}_{n} \right]= \hat{X}_{n}\\
\left[ \hat{X}_{i}, \hat{T}_{j} \right] = - a_{ij} \hat{X}_{i}
\end{array}
\end{equation}
and the quantum Serre relations described by the expressions:
 \begin{equation}
 \begin{array}{c}
  \sum_{0\le t \le 1-a_{ik}} (-1)^{t}
 \left[  \begin{array}{clcr}
 &1-a_{ik}\\ &t \end{array}  \right]_{Q}\hat{X}_{i}^{1-a_{ik}-t} \hat{X}_{k} \hat{X}_{i}^{t}= 0\\
 \sum_{0\le t \le 1-a_{ij}} (-1)^{t} \left[  \begin{array}{clcr}
 &1-a_{ij} \\&t \end{array}  \right]_{Q}\hat{Y}_{i}^{1-a_{ij}-t} \hat{Y}_{j} \hat{Y}_{i}^{t} =0
 \end{array}
  \end{equation}
A possible realization of the generators of $ \tilde{P}_{Q}(n)$ in terms of the $Q$-deformed oscillators $ \{ \Phi_{i}^-,\Phi_{i}^+$, $M_{\Phi_i} \}$  and $\{A^{+}_{i},A^{-}_{i},N_{A_{i}} \}$ with $( 1 \leq i \leq n )$ is given by:
\begin{equation}
 \begin{array}{c}

 \hat{X}_{i}= A^{+}_{i}A^{-}_{i+1}Q^{\frac{(M_{\Phi_{i}}-M_{\Phi_{i+1}})}{2}} + \Phi^{+}_{i} \Phi^{-}_{i+1}Q^{\frac{-(N_{A_{i}}-N_{A_{i+1}})}{2}}\\\\

\hat{X}_{n}= A^{+}_{n}\Phi^{+}_{n}Q^{(\frac{1}{2} \sum_{i=1}^{n-1} N_{A_{i}} - \frac{1}{2} \sum_{i=1}^{n-1} M_{\Phi_{i}})} \\\\

\hat{Y}_{i}= A^{+}_{i+l}A^{-}_{i}Q^{\frac{(M_{A_{i}}-M_{A_{i+1}})}{2}} + \Phi^{+}_{i+1} \Phi^{-}_{i}Q^{\frac{-(N_{\Phi_{i}}-N_{\Phi_{i+1}})}{2}}\\\\

\hat{T}_{i}= N_{A_{i}} - N_{A_{i+1}} + M_{\Phi_{i}} - M_{\Phi_{i+1}}\\\\

\hat{D} =\frac{1}{2} \sum_{i=1}^{n} N_{A_{i}} + \frac{1}{2} \sum_{i=1}^{n} M_{\Phi_{i}}\\\\

\end{array}
 \end{equation}
where $A^{+}_{i}$ and $A^{-}_{i}$ are the $Q$-deformed bosonic algebra operators and $\Phi^{+}_{i}$and $\Phi^{-}_{i}$ are the $Q$-deformed
fermionic ones.\\
Using  the fact that in $Q\to q ~limit$, the
$Q$-deformed boson algebras $\{A^{+}_{i},A^{-}_{i},N_{A_{i}} \}$ reproduces  classical bosons algebras $%
\{a^{+}_{i},a^{-}_{i},N_{a_{i}}\}$ and a q-deformed k-fermions generators $\{\chi^{+}_{i},\chi^{-}_{i},N_{\chi_{i}}\},$ and in this
$limit$ the $Q$-fermions become $q$-fermions which are object equivalent to
classical ones $\{\phi^{+}_{i},\phi^{-}_{i},M_{\phi_{i}}\}$, then from the classical bosons $\{a_{i}^{-}, a_{i}^{+}, N_{a_{i}} \}$ and classical fermions $\{\phi_{i}^{-}, \phi_{i}^{+}, M_{\phi_{i}} \}$, one can construct
  the classical extended strange superalgebra $\tilde{P}(n)$ :
 \begin{equation}
 \begin{array}{c}
 X_{i}= a^{+}_{i}a^{-}_{i+1} + \phi^{+}_{i} \phi^{-}_{i+1}\\
 X_{n}= a^{+}_{n}\phi^{+}_{n} \\
 Y_{i}= a^{+}_{i+l}a^{-}_{i} + \phi^{+}_{i+1} \phi^{-}_{i}\\
 T_{i}= N_{a_{i}} - N_{a_{i+1}} + M_{\phi_{i}} - M_{\phi_{i+1}}\\
 D =\frac{1}{2} \sum_{i=1}^{n} N_{a_{i}} + \frac{1}{2} \sum_{i=1}^{n} M_{\phi_{i}}\\
 \end{array}
 \end{equation}
  From the remaining operators $\{\chi_{i}^{-}, \chi_{i}^{-}, N_{\chi_{i}} \}$ we realize the following :
 \begin{equation}
 \begin{array}{c}
 E_{i}= \chi_{i}^{-} \chi_{i+1}^{+},~~~ 1 \le i \le n-1\\
 F_{i}= \chi_{i}^{+} \chi_{i+1}^{-},~~~ 1 \le i \le n-1\\
K_{i}= q^{-N_{\chi_{i}}+ N_{\chi_{i+1}}},~~~ 1\le i \le n-1\\
 \end{array}
 \end{equation}
  which generates the $q$-deformed algebra $U_{q}(sl(n))$. It is easy to show that $U_{q}(sl(n))$ and $\tilde{P}(n)$ are mutually commutative. As results, we obtain the following decomposition of the quantum strange superalgebra $\tilde{P}_{Q}(n)$ in the $Q \to q~ limit$
 \begin{equation}\label{edec}
   lim_{Q \to q}\tilde{P}_{Q}(n) = U_{q}(sl(n)) \otimes \tilde{P}(n).
\end{equation}

\section{Conclusion}
It is important to note that we have established this decomposition of the quantum strange superalgebra $\tilde{P}_{Q}(n)$ only for a particular realization, i.e, the $Q$-oscillator realization and although the quantum extended strange superalgebra $\tilde{P}_{Q}(n)$ does not have direct product form, we establish, for this realization and the corresponding highest weight representations the decomposition of $\tilde{P}_{Q}(n)$ into the direct product of undeformed $\tilde{P}(n)$ and $U_{q}(sl(n))$ (the naive version of $U_{Q}(sl(n))$  at $Q= q $ obtained by simply setting $Q= q $). The labels of the highest weight representations of the quantum strange superalgebra $\tilde{P}_{Q}(n)$ and the choice of the basis in which the decomposition (\ref{edec}) is  clearly manifested will be investigated elsewhere.

\newpage\

\end{document}